\magnification 1200
\input amssym.def
\input amssym.tex
\parindent = 40 pt
\parskip = 12 pt
\font \heading = cmbx10 at 12 true pt
 at 22 true pt
 at 19 true pt
 at 7 true pt
\def \R{{\bf R}}

\centerline{\heading An Elementary Proof of the Meromorphy}
\centerline{\heading of $\int f^z$ for real-analytic $f$}
\rm
\line{}
\line{}
\centerline{\heading Michael Greenblatt}
\centerline{mg62@buffalo.edu}
\line{}
\centerline{January 26, 2007}
\baselineskip = 12 pt
\font \heading = cmbx10 at 14 true pt
\line{}
\line{}

Suppose $U$ is a domain in $\R^n$, $f$ is a real-analytic function on $U$ and $M$ is an open
semianalytic subset of $U$ with $cl(M) \subset U$. For $\phi \in C_c(U)$, we define
$$F(z) = \int_M (f(x))^z \phi(x) \, dx \eqno (1.1)$$
To be clear what $(1.1)$ means when $f$ takes on negative values, we fix some branch of the logarithm 
on the negative real axis and use it
to define $f(x)^z$ when $f(x) < 0$. If ${\rm Re}(z) > 0$, the integrand in $(1.1)$ is
integrable and standard methods show that $F$ is analytic in a neighborhood of $z$ with 
$$F'(z) = \int_M \log(f(x))(f(x))^z \phi(x) \, dx \eqno (1.2)$$
Here $\log(f(x))$ denotes the branch of the logarithm chosen above. 
It is natural to try to extend $F(z)$ to a meromorphic function of $z$ on all of ${\bf C}$ by doing
appropriate integrations by parts, integrating $f(x)^z$ in some way. However, the zero set of an 
arbitrary real-analytic function can be quite complicated, so carrying this out may be rather involved.
After Hironaka proved his monumental [H1] [H2], Gelfand suggested that by using these results one might 
be able to do the requisite analysis by reducing to the case where $f(x)$ is a monomial. This was 
done by Bernstein and Gelfand [BGe] when $f(x)= |P(x)|$ for a polynomial $P$ and when $M = U
= \R^n$. In addition, Atiyah [A] did it for general semianalytic $M$ and general nonnegative
real-analytic $f(x)$. Later, Bernstein [B] found an algebraic proof of the results in [BGe] not 
using resolution of singularities, by virtue of Bernstein-Sato polynomials, which allow one to 
integrate by parts directly in the original integral.

The purpose of this paper is to show that the elementary resolution of singularities algorithm of 
[G1] suffices to prove such results. In fact, we prove a very slightly more general result than 
that of [A] by 
removing the restriction that $f$ is nonnegative (the version of resolution of singularities used in 
[A] actually does use that $f \geq 0$ on all of $U$, so one can't simply resolve $f^2$ or add a 
condition like $f < 0$ to the definition of $M$). The main result of this paper is:

\noindent {\bf Theorem 1:} The function $F(z)$ of $(1.1)$ extends to a meromorphic function on
{\bf C}. If $K$ denotes any compact set containing $supp(\phi)$, then poles of the resulting 
function must be at a point of the form $-{r \over N}$,
where $N$ is a fixed positive integer depending on $f$, $M$, and $K$, where $r$ is a 
positive integer. The order of any pole is at most the dimension $n$. 

\noindent To prove Theorem 1, we use the following consequence of the Main Theorem of [G]:

\noindent {\bf Resolution of Singularities Theorem:} Suppose $f$ is real analytic on a neighborhood of the origin. Then there is
an neighborhood $V$ of the origin such that if $\phi(x) \in C_c(V)$ is nonnegative with $\phi(0) > 0$,
then $\phi(x)$ can be written as $\phi(x) = \sum_{i=1}^p \phi_i(x)$, each $\phi_i(x)$ nonnegative,
such that the following hold. Let $D_i = \{x: \phi_i(x) > 0\}$. There is a real-analytic diffeormorphism
$\Psi_i$ from an open bounded $D_i'$ to $D_i$ such that on a neighborhood of $cl(D_i')$, 
$f \circ \Psi_i(x) = d_i(x)m_i(x)$, $m_i(x)$ a monomial and $d_i(x)$ nonvanishing.  

\noindent One can resolve several functions simultaneously in this fashion by resolving their product. 

\noindent The Main Theorem of [G] also stipulates that each $\phi_i \circ \Psi_i(x)$ is a "quasibump function":

\noindent {\bf Definition:} Let $E = \{x : x_i > 0$ for all $i\}$. If $h(x)$ is 
a bounded, nonnegative, compactly supported function on $E$, we say $h(x)$ 
is a {\it quasibump function} if $h(x)$ is of the following form:
$$h(x) = a(x) \prod_{k=1}^l b_k (c_k(x) {p_k(x) \over q_k(x)}) \eqno (1.2)$$
Here $p_k(x), q_k(x)$ are monomials, $a(x) \in C^{\infty}(cl(E))$, the $c_k(x)$ are  
nonvanishing real-analytic functions defined
on a neighborhood of $supp(h)$, and $b_k(x)$ are
functions in $C^{\infty}(\R)$ such that there are  $c_1 > c_0 > 0$ with each $b_k(x)
= 1$ for $x < c_0$ and $b_k(x) = 0$ for $x > c_1$. 

The reason this explicit form for $h(x)$ is useful for our purposes is that we
will be doing some integrations by parts in integrals with $\phi_i \circ \Psi_i(x)$ appearing in the 
integrands, and derivatives will be landing on such $\phi_i\circ \Psi_i(x)$. Since we
know the $\phi_i\circ \Psi_i(x)$ are all quasibump functions, we can get explicit estimates on the size of the 
derivatives. In addition, since the $b_k(x)$ are constant on $x < c_0$ and on $x > c_1$, the 
support of a derivative of $\phi_i\circ \Psi_i(x)$ will be substantially smaller than that of 
$\phi_i \circ \Psi_i$. After an
appropriate coordinate change, this will effectively allow us to reduce the dimension of the 
problem and induct on the dimension $n$.

\noindent {\bf Proof of Theorem 1:} We can assume that $M$ is of the form $\{x : g_k(x) > 0$ for 
$k = 1,...,p \}$, each $g_k$ real-analytic, since up to a set 
of measure zero every semianalytic set can be written  as the 
finite union of sets of this form. At each point $x$ in $supp(\phi) \cap cl(M)$ one can find a
neighborhood $N_x$ of $x$ such that the resolution of singularities theorem above applies to 
$f, g_1,...,g_p$ simultaneously on $N_x$.
Using a partition of unity, on $supp(\phi) \cap cl(M)$ one can write $\phi = \sum_j
\alpha_j(x)$, where each $\alpha_j$ is in $C_c(N_x)$ for some $x$. Hence it suffices to prove Theorem
1 for an arbitrary $\alpha_j$ in place of $\phi$; adding the results will give Theorem 1 for
$\phi$. Hence without loss of generality, we assume $\phi$ is one of these $\alpha_j$. Thus 
we may apply the resolution of singularities theorem to $f, g_1,...,g_p$, and we write $\phi = 
\sum \phi_i$ accordingly. Define
$$F_i(z) = \int_M (f(x))^z \phi_i(x) \, dx \eqno (1.3)$$
It suffices to show each $F_i(z)$ satsifies the conclusions of Theorem 1. We change 
coordinates in $(1.3)$ to the blown up coordinates, obtaining
$$F_i(z) = \int_E (f\circ \Psi_i(x))^z (\chi_M \circ \Psi_i(x)) (\phi_i \circ \Psi_i(x))\,det\, \Psi_i(x) 
\, dx \eqno (1.4)$$
In the new coordinates, we have $f \circ \Psi_i(x) = d_i(x)m_i(x)$, where $d_i(x)$ is nonvanishing
on a neighborhood of $supp(\phi_i \circ \Psi_i)$ and $m_i(x)$ is some monomial $\prod_{j=1}^n x_j^{a_{ij}}$. Hence we 
can rewrite $(1.4)$ as 
$$F_i(z) = \int_E \prod_{j=1}^n x_j^{a_{ij}z} d_i(x)^z (\chi_M \circ \Psi_i(x))( \phi_i \circ \Psi_i(x))
\,det \,\Psi_i(x) \, dx \eqno (1.5)$$
In the case that $d_i(x)$ assumes negative values, one defines $d_i(x)^z$ using the branch of 
the logarithm one used to define $f(x)^z$ in the original coordinates. Since each $g_k \circ 
\Psi_i(x)$ is also of the form $d(x)m(x)$, either each $g_k$ is positive throughout the 
domain of integration of $(1.5)$, or there is at least one $k$ for which $g_k$ is negative
throughout the domain. In the first case, $\chi_M \circ \Psi_i(x)$ is always 1, in the second 
case it is always zero. In the latter case $F_i(z) = 0$ and there is nothing to prove, so we assume
we are in the first case and write
$$F_i(z) = \int_E \prod_{j=1}^n x_j^{a_{ij}z} d_i(x)^z  (\phi_i \circ \Psi_i(x)) \,det\,\Psi_i(x) 
\, dx $$
We use the explicit form $(1.2)$ of the quasibump function $\phi_i \circ \Psi_i(x)$ and this becomes
$$F_i(z) = \int_E \prod_{j=1}^n x_j^{a_{ij}z} d_i(x)^z  a(x) \prod_{k=1}^l b_k (c_k(x) {p_k(x) 
\over q_k(x)})\,det\,\Psi_i(x)\,dx$$
We combine the smooth $a(x)$ and $det\,\Psi_i(x)$ factors by letting $A_i(x) = a(x)\,det\,\Psi_i(x)$, and
the above becomes
$$F_i(z) = \int_E \prod_{j=1}^n x_j^{a_{ij}z} d_i(x)^z  A_i(x) \prod_{k=1}^l b_k (c_k(x) {p_k(x) 
\over q_k(x)})\,dx\eqno (1.5')$$
The idea behind the analysis of $(1.5')$ is quite simple. One wishes to repeatedly integrate by 
parts, integrating first $x_1^{a_{i1}z}$, then $x_2^{a_{i2}z}$, going up to $x_n^{a_{in}z}$. Then
one integrates $x_1^{a_{i1}z + 1}$, cycles through the $x_j$ again, and repeats ad nauseum. (One may
skip any $x_j$ for which $a_{ij} = 0$). Since the
exponents of the $x_j$ will increase each time, the integral will be analytic over a larger
and larger $z$-domain as the integrations by parts proceed. One obtains poles since one gets
factors of ${1 \over a_{ij}z + k}$ showing up with each integration by parts. Each pole should
have order at most $n$, since a given factor appears at most once per variable. 

To make these heuristics work, one has to ensure that the derivatives landing on the  
$b_k (c_k(x) {p_k(x) \over q_k(x)})$ factors don't mess things up. In [A] or [BGe] such issues 
don't arise since they use the stronger result of Hironaka which
doesn't require one to subdivide a neighborhood of the origin into different parts each 
having a different set of coordinate changes; instead there is one sequence $\Psi$ of blow ups and 
the resulting $\phi \circ \Psi$ is smooth.

\noindent Theorem 1 will follow from the following lemma:

\noindent {\bf Lemma:}  $F_i(z)$ extends to a 
meromorphic function on ${\bf C}$. There is a positive integer $N$, depending on 
the $a_{ij}$ and the various monomials $p_k(x)$ and $q_k(x)$, such that each 
pole of $F_i(z)$ is at $-{r \over N}$ for some nonnegative integer $r$. Let $\eta > 0$ such that 
the integrand of $(1.5')$ is supported on $(0,\eta)^n$. Then for any $l$ and  
each compact subset $K$ of $\{z: z > -{l + 1 \over N}\}$, 
$\sup_K\vert \prod_{r=1}^l(z + {r \over N})^n F_i(z)\vert $ can be bounded in terms of $K$,
$\eta$, the $a_{ij}$, the monomials $p_k(x)$ and $q_k(x)$, and the $C^m$ norms of the
$d_i(x)$, $A_i(x)$, $b_k(x)$, and $c_k(x)$. Here $m$ is some sufficiently large natural number. 

\noindent {\bf Proof:} We proceed by induction on $n$. We do the $n=1$ case at the same time
as the $n > 1$. So we assume that either $n = 1$ or that $n > 1$ and we know the result for
$n - 1$.  We perform an integration by parts in $x_1$ in
$(1.5')$, turning the $x_1^{a_{i1}z}$ into
an ${x_1^{a_{i1}z + 1} \over {a_{i1}z + 1}}$. If the derivative lands on $d_i(x)^z$, we obtain
another smooth function which does not interfere with future integrations by parts with respect to
$x_2$, $x_3$, etc as described in the heuristics above. A similar situation
occurs if the derivative lands on $A_i(x)$. Things are more complicated when the derivative
lands on one of the 
$b_k (c_k(x) {p_k(x) \over q_k(x)})$. Let $x_1^m$ denote the power of $x_1$ appearing in 
${p_k(x) \over q_k(x)}$ (hence $m$ could be positive, negative, or zero). Then we have
$$\partial_{x_1}[b_k (c_k(x) {p_k(x) \over q_k(x)})] = {1 \over x_1}(x_1 {\partial c_k \over 
\partial x_1}(x) + mc_k(x)) ({p_k(x) \over q_k(x)})b_k'(c_k(x) {p_k(x) \over q_k(x)}) \eqno (1.6)$$
Thus if we write $B_k(x) = xb_k'(x)$, we have
$$ \partial_{x_1}[b_k (c_k(x) {p_k(x) \over q_k(x)})] = {1 \over x_1}({x_1 {\partial c_k \over 
\partial x_1}(x) + m \over c_k(x)})B_k (c_k(x) {p_k(x) \over q_k(x)}) \eqno (1.7)$$
We consider $({x_1 {\partial c_k \over \partial x_1}(x) + m \over c_k(x)}) a(x)$ as
a smooth factor $s_k(x)$, and thus we have   
$$\partial_{x_1}[b_k (c_k(x) {p_k(x) \over q_k(x)})] = {1 \over x_1}s_k(x) B_k (c_k(x) 
{p_k(x) \over q_k(x)})\eqno (1.8)$$ 
If we use the notation $A_{ik}(x) = A_i(x)s_k(x)$, the integral corresponding to this term is given by
$${1 \over a_{i1}z + 1} \int_E \prod_j x_j^{a_{ij}z} d_i(x)^z [A_{ik}(x)B_k (c_k(x) 
{p_k(x) \over q_k(x)}) \prod_{K \neq k} b_K (c_K(x) {p_K(x) \over q_K(x)})]  \,
dx \eqno (1.9)$$
At first glance, $(1.9)$ might appear to be little improved over $(1.5')$, since the exponents
of the $a_{ij}z$ appearing are unchanged.  However, there is a key difference. Namely, instead of
having the quasibump function $\phi_i \circ \Psi_i(x)$ in the integrand, we have the bracketed
expression in $(1.9)$.
Because $b_k(x)$ is constant for $0 < x < x_0$ and
$x > x_1$, $B_k(x) = xb_k'(x)$ is supported on $[x_0,x_1]$. This means that the factor 
$B_k (c_k(x) {p_k(x) \over q_k(x)})$ in the integrand is supported on the wedge $C_1q_k(x)
\leq p_k(x) \leq C_2q_k(x)$ for some constants $C_1$ and $C_2$. After doing appropriate coordinate
changes, we will be able to exploit this fact to reduce the problem to the $n-1$ dimensional
case. We break into three
cases of increasing order of difficulty. We can assume $p_k(x)$ and $q_k(x)$ have no common 
factors. 

\noindent Case 1): Either $p_k(x)$ or $q_k(x)$ is constant, and each $x_j$ appears  to a 
positive power in whichever of $p_k(x)$ or $q_k(x)$ is nonconstant. Note that whenever
$n = 1$ we are in case 1. Replacing $B_k(x)$ by $B_k(1/x)$ if
necessary, we may assume $p_k(x)$ is constant. 
Because $B_k(x)$ is zero for $x > x_0 > 0$ for some $x_0$ and the integrand of 
$(1.9)$ is compactly supported, there is some constant $C$ such that each $x_j > C$
when the integrand of $(1.9)$ is nonzero. As a result, one can integrate by parts in $(1.9)$ as
many times as one likes; the integral is in fact an entire function. 

\noindent Case 2): Either $p_k(x)$ or $q_k(x)$ is constant, but there is some $x_j$ not
appearing in the nonconstant function of $p_k(x)$ or $q_k(x)$. Like before we may replace 
$B_k(x)$ by $B_k(1/x)$ if necessary and assume
$p_k(x)$ is constant. Let $J$ be the set of $j$ for which $x_j$ appears in $q_k(x)$. Then 
in the integrand of $(1.9)$, $x_j > C$ for all $j \in J$. Thus if in the integrand of $(1.9)$ we
freeze each $x_j$ at a constant for $j \in J$, then the integrand becomes that of an expression
$(1.5')$ corresponding to the $n - |J|$ dimensional case; we treat $A_{ik}(x)B_k (c_k(x) 
{p_k(x) \over q_k(x)})$ like a single smooth factor $A_i(x)$.
Hence by induction hypothesis, the integral in these
$n - |J|$ variables is a meromorphic function satisfying the conclusions of the lemma in 
dimension $n - |J|$. By the uniform bounds given by 
the lemma, we conclude that $(1.9)$ satisfies the conclusions of the lemma as well, and we are
done with case 2.

\noindent Case 3): Both $p_k(x)$ and $q_k(x)$ are nonconstant. In this case we will have to
break up $(1.9)$ into several pieces. Some coordinate changes are done on each piece to reduce
it to Cases 1 or 2. We first do a coordinate change $(x_1,...,x_m) \rightarrow (x_1^{M_1},...,
x_n^{M_n})$ so that each $x_j$ appearing in either $p_k(x)$ or $q_k(x)$ appears to the
same power. We still get an expression of the form $(1.9)$, after incorporating
the determinant of this coordinate change into the $A_{ik}(x)$ factor.
Suppose $x_l$ appears in $p_k(x)$ and $x_m$ appears in $q_k(x)$. Let $\alpha(x)
\in C^{\infty}(0,\infty)$ be nonnegative such that $\alpha(y) + \alpha(1/y) = 1$ for all $y$,
and such that $\alpha(y)$ is supported on $y < C$ for some $C$. In particular,
in the integrand of $(1.9)$ we have $\alpha({x_l \over x_m}) + \alpha({x_m \over x_l}) = 1$, and
we correspondingly write the integral as  $I_1(z) + I_2(z)$, where 
$$I_1(z) = \int_E \prod_j x_j^{a_{ij}z} d_i(x)^z [A_{ik}(x)B_k (c_k(x) 
{p_k(x) \over q_k(x)}) \prod_{K \neq k} b_K (c_K(x) {p_K(x) \over q_K(x)})]\alpha({x_l \over x_m})
 \, dx \eqno (1.10)$$
$$I_2(z) = \int_E \prod_j x_j^{a_{ij}z} d_i(x)^z [A_{ik}(x)B_k (c_k(x) 
{p_k(x) \over q_k(x)}) \prod_{K \neq k} b_K (c_K(x) {p_K(x) \over q_K(x)})]\alpha({x_m \over x_l})
 \, dx \eqno (1.11)$$
In $(1.10)$ we do the variable change turning what was $x_l$ into $x_l x_m$, and in $(1.11)$ we
do the variable change turning what was $x_m$ into $x_m x_l$. The resulting integrals are still
of the form $(1.9)$. In addition, in the factor ${p_k(x) \over q_k(x)}$ of $(1.10)$ the factor
$x_m$ dissapears, while in $(1.11)$ the factor $x_l$ dissappears. 

If we iterate the above in $(1.10)$ and $(1.11)$, splitting into more and more terms, then 
eventually enough $x_j$'s will have dissappeared in ${p_k(x) \over q_k(x)}$ that we are in either case
1 or case 2. Thus $(1.9)$ is the sum of finitely many terms that fall under case 1 or case 2, and
thus we have the lemma in case 3 as well. This completes the proof of the lemma. We are also
done with the proof of Theorem 1; as in the earlier heuristics we integrate by parts with respect to 
$x_2$, $x_3$, etc ad infinitum; one deals with these integrations by parts the way we dealt with
the $x_1$ integration by parts above. 

\noindent {\bf References:}

\noindent [A] M. Atiyah, {\it Resolution of singularities and division of distributions},
Comm. Pure Appl. Math. {\bf 23} (1970), 145-150.

\noindent [BGe] I. N Bernstein and S. I. Gelfand, {\it Meromorphy of the function ${\rm P}^
\lambda$}, Funkcional. Anal. i Prilo\v zen. {\bf 3} (1969), no. 1, 84-85. 

\noindent [B] I. N. Bernstein, {\it Analytic continuation of generalized functions with respect 
to a parameter}, Funkcional. Anal. i Prilo\v zen. {\bf 6} (1972), no. 4, 26-40.

\noindent [G1] M. Greenblatt, {\it A Coordinate-Dependent Local Resolution of Singularities and 
Applications}, preprint.

\noindent [H1] H. Hironaka, {\it Resolution of singularities of an algebraic 
variety over a field of characteristic zero I},  Ann. of Math. (2) {\bf 79}
(1964), 109-203;

\noindent [H2] H. Hironaka, {\it Resolution of singularities of an algebraic 
variety over a field of characteristic zero II},  Ann. of Math. (2) {\bf 79}
(1964), 205-326.

\end